\documentstyle[12pt]{amsart}

\theoremstyle{plain}
 \newtheorem{thm}{Theorem}[section]
 \newtheorem{prop}[thm]{Proposition}
 \newtheorem{lemma}[thm]{Lemma}
 \newtheorem{cor}[thm]{Corollary}
\theoremstyle{definition}

\theoremstyle{remark}

\begin{document}
\title{Quantized  W-algebra of   ${\frak sl}(2,1)$ : \\
a construction from the quantization of 
screening operators}  
\author{Jintai Ding}
\address{Jintai Ding, RIMS, Kyoto University}
\author{Boris Feigin}
\address{Landau Institute for Theoretical
Physics, Moscow;  Russia and RIMS, Kyoto University}
\maketitle
\begin{abstract}
Starting from bosonization,  we study the operator that 
commute or commute up-to a total difference 
  with of any  quantized  
screen operator of  a free field. We show that  if there exists 
a  operator in the form of a sum of  two vertex operators
which has  the simplest  correlation functions  with the 
quantized screen operator, 
namely a function with one pole and one zero, then,  
the screen operator and this operator are  uniquely determined, and
this operator  is the quantized virasoro algebra.
For the case when the screen is a fermion, there are a family of this kind
of operator, which give new algebraic structures. Similarly we 
 study  the case of two quantized  
screen operator, which uniquely gives us the quantized W-algebra 
corresponding to $sl(3)$ for the generic case, and a new algebra, which 
is  a quantized W-algebra  corresponding to ${\frak sl}(2,1)$, for the case
 that 
one of the two  screening operators  is or both are fermions. 
\end{abstract}

\section{Introduction.}

In \cite{FR} N. Reshetikhin and E. Frenkel 
introduced new Poisson
algebras $W_q({\frak g})$, which are $q$--deformations of the classical
$W$--algebras. 
In \cite{SKAO} J. Shiraishi, H.Kubo, H.Awata, and S. Odake quantized
the formulation in \cite{FR}.  They constructed a
non-commutative algebra depending on two parameters $q$ and $p$, such that
when $q=p$ it becomes commutative, and is isomorphic to the Poisson algebra
$W_q({\frak sl }(2))$, which denote by  $W_{q,p}({\frak sl}(2))$. Shiraishi,
e.a., constructed a free field realization of , $W_q({\frak sl }(2))$ i.e. a
homomorphism into a Heisenberg algebra. 
They also constructed the screening currents,
i.e. operators acting on the Fock representations of the Heisenberg algebra.
which commute with the action of $W_q({\frak sl }(2))$ up to a total
difference. 
In \cite{FF}, the results of Shiraishi,
e.a., were further 
generalized to 
to the case of the $q$-deformed $W$-algebras. They 
constructed  an algebra $W_{q,p}({\frak sl }(N))$ 
depending on $q$ and $p$, such that
when $q=p$ it becomes isomorphic to the $q$--deformed classical
$W$--algebra $W_q({\frak sl }(N))$ from \cite{FR}. They also construct
 a free field realization of
$W_{q,p}({\frak sl }(N))$, 
which is a deformation of the free field realization
from \cite{FR}, and the screening currents. 
They also construct the 
screening currents $S^\pm_i(z)$ satisfying certain 
difference equations related to the basic generator of $W_{q,p}(
{\frak sl }(N))$. 

In all the works above, the screening operators, 
which are quantization of the classical screening operators, 
 appear as a by-product 
of the quantization of  W-algebras. These 
 screening operators are 
defined as  operators who
either commute with the quantized W-algebras or 
commute  up-to a total difference.

In this paper, we will start from the opposite  direction of 
such an  approach. Namely, we 
will start from the quantization of 
screening  operators. This approach is 
based on the idea of bosonization, namely we start everything 
from the Heisenberg algebras. 

First we consider the case of 
one generic screen operator (not a fermion)
$S^+_1(z)$, which is given as 
any kind of   quantization  of the classical screen operator 
$S(z)$. Here both operator are bosonized, namely expressed by 
generator of a  Heisenberg algebra.  
Surely the screening operator $S(z)$ satisfies: 
 $$S(z)S(w)= (z-w)^{2\beta} :S(z)S(w):,  $$
and 
 $$S^+_1(z)S^+_1(w)= (z)^{2\beta}f(z,w) :S^+_1(z)S^+_1(w) , $$
where $\beta \in \Bbb C$. 
By quantization, we mean that the operator $S^+_1(z)$ depends on a
parameter q and it degenerates into the classical 
screening operator when q goes to 1 or 
$$\lim_{q\rightarrow 1} z^{2\beta}f(z,w)=(z-w)^{2\beta} . $$
   Then we  try to construct an  operator $l(z)=\Lambda_1(z)+\Lambda_2(z)$, 
a  sum of two  vertex operators, that will commute with 
the action of the screen operator up-to a total 
difference. It turns out that 
this uniquely determines the 
quantization  and  enable us to   recover the quantized virasoro algebra 
mentioned above in an unique way,  if we assume that the correlation 
functions between $S^+_1(z)$ and $\Lambda_i(w)$ has one pole and one zero
and the ration of the two pole is our quantization parameter $q$ .  
However, when the screen operator is a fermion, the situation becomes very 
different. For this case, we will obtain a family of 
operators that commute with the screening operator up-to a
 total difference, which give  completely  new algebraic structures.
Similarly, we consider the case of two 
generic screen operators  (neither of them  a fermion), which are given as 
any kind of   q-deformation of the classical screen operators.
Then we try to construct an operator in the form of 
the sum  of three vertex operator, such that this operator commutes 
with the two screening operators. This enable us to uniquely derive 
the quantized W-algebra for $sl(3)$ given by Feigin and Frenkel 
\cite{FF}.    

Then we apply the same method to the case of two screening operators 
such that one of the screening operators is a fermion. In this case, 
this leads us to derive an unique operator such that it commutes
with  the two screening operator up-to 
a total difference. We will call the algebraic 
structure generated by this operator the quantized W-algebra of 
${\frak sl}(2,1)$. Finally we apply this method to the case that 
two screening operators are both fermions. This also leads us to 
derive an unique operator that  generates the same algebraic 
structure, namely the  quantized W-algebra of 
${\frak sl}(2,1)$. 

This paper is organized in the way that  each section is devoted to study 
a case mention above and in the same order. The last sections 
is devoted to discussions.


\section{ The commutant of a screen operator and the 
q-deformed virasoro algebra}    

In this section, we will start from one 
quantized screen operator.  From the point view of quantization, 
this screen operator should degenerate into the classical screen 
operator, when the deformation parameter $q$ goes to $1$. With 
such a screen operator, we will try to find current operators, 
which commute with the integral of this screen operator. 
Analog  to the idea that differential operators in general are deformed into 
difference operators, we expect that 
a desired operator, which should be a kind of q-deformed virasoro algebra, 
 is in the  form of a kind of difference.

Let us first introduce the
Heisenberg  algebra 
$H_{q,p}(1)$ be the Heisenberg algebra with generators $a_1[n],
n\in{\Bbb Z}$, and relations
\begin{equation}    
\left[ a_1[n],a_1[m] \right] = \frac{1}{n}
\delta_{n,-m},
\end{equation}
which is defined on the filed of the rational functions of 
 $p=q^x$ and $q$. Here $x$ and $p$ 
are two generic parameters  and $|q|<1$

For each weight $\mu$ of the Cartan subalgebra of ${\frak sl}(2)$, 
let $\pi_\mu$ be
the Fock representation of $H_{q,p}(1)$ generated by a vector
$v_\mu$, such that $a_1[n] v_\mu = 0, n>0$, and $a_1[0] v_\mu =
\mu(\alpha_1^\vee) v_\mu$, 
where $\alpha_1^\vee$ is the $i$th coroot of ${\frak sl}(2)$.

Introduce operators $Q_1$, which satisfy commutation
relations $[a_1[n],Q_1] = 2\beta \delta_{n,0}$. The operators
$e^{Q_1}$ act from $\pi_\mu$ to $\pi_{\mu+\beta\alpha_1}$.

A classical screening operator $S(z)$ is defined as
\begin{align}    
S(z) & = e^{Q_1} z^{a^+_1[0]} :\exp \left( \sum_{m\neq 0} s(m)
a_1[m]z^{-m} \right):, 
\end{align}
where $m\neq 0$ and  
\begin{align} 
S(z)S(w)= (z-w)^{2\beta}:S(z)S(w):. 
\end{align} 
 
Now we can define a quantized screening currents as the generating function
\begin{align}    
S^+_1(z) & = e^{Q_1} z^{s^+_1[0]} :\exp \left( \sum_{m\neq 0} s^+_1(m)
a_1[m]z^{-m} \right):, 
\end{align}
where $s^+_1[m]$ are in $\Bbb C[p,q]$ 
for $m\neq 0$ and $s^+_1[0] = a_1[0]$. 
For this screen operator we would impose the following condition that 
the limit of this operator when $q$ goes to one degenerate into the 
classical screen operators.

For the undeformed case, the virasoro algebra defined on the Fock space
give us the operators that commute with the action of the operator 
$\int S^+_1(z)dz/z$. Here we will impose the 
assumption  that once we consider 
such an integration,  the integration contour is around the point 0 
and the  screening operator $S^+_1(z)$ acts on the
space   $\pi_\mu$ such that  this integration 
is single-valued. 

We, then,  would try  to find a 
construction that gives us the operators that will commute with the 
action of this quantized screen operator $\int S^+_1(z)dz/z$. 
The simplest try would be that 
case that they should be a sum of vertex operators. 

For the simplest case, we assume that it is sum of two 
vertex operators. Let us define the  operator as 
$$l_1(z)= \Lambda_1(z)+\Lambda_2(z),$$  where 
$\Lambda_i(z)$ as the generating function: 
\begin{align}    
\Lambda_i(z) & =g^{i-1}
 :\exp \left( \sum_{} \lambda_i(m)
a_1[m]z^{-m} \right):.
\end{align}
Here $\lambda_i[m]$ are in $\Bbb C[p,q]$  for $i=1,2$. 

The  commutation relations between $S^+(z)$ and 
$\Lambda_i(w)$ are basically decided by the correlation functions 
of the product of this two operators. The simplest case is that 
the correlation function of 
 the operators  has  only one pole. Similarly,
we assume that it has at most one zero. 
Clearly, it naturally leads to that we must 
have the condition that the two 
products  $\Lambda_1(z) S^+_1(w)$
and  $S^+_1(w) \Lambda_1(z)$ have the same correlation functions. 
Then we have 
$$\Lambda_1(z) S^+_1(w) =
A\frac{(z-wp_1)}{(z-wp_2)} :\Lambda_1(z)
S^+_1(w):, \quad \quad |z|\gg|w|,$$ $$S^+_1(w) \Lambda_1(z) =
A\frac{(z-p_1w)}{(z-p_2w)} :\Lambda_1(z)
S^+_1(w):, \quad \quad |w|\gg|z|.$$ 
This implies:

\begin{lemma}
$$Ap_1/p_2=1.$$
\end{lemma}
Therefore $p_1$ can not be 0. 

The correlation functions also give us 
$$\Lambda_1(z) S^+_1(w)= A
\exp \left (\Sigma_{n>0}\lambda_1(n)s^+_1(-n)(w/z)^n \right) :\Lambda_1(z)
S^+_1(w):$$
$$ = A
\exp \left (\Sigma_{n>0}(w/z)^n(-p_1^n+p_2^n) \right) :\Lambda_1(z)
S^+_1(w):.$$
Thus we have 
$$\lambda_1(n)s^+_1(-n)= (-p_1^n+p_2^n).$$
Similarly we have 
$$\lambda_1(-n)s^+_1(n)= (-p_1^{-n}+p_2^{-n}). $$

\begin{lemma}
If the operator $l_1(z)$ commute with the 
the operator   $\int S^+_1(z)dz/z$, then 
 $l_1(z)$  commute with $S^+_(w)$ up-to a total difference, 
the correlation functions of the products 
$\Lambda_1(z) S^+_1(w)$
and  $S^+_1(w) \Lambda_1(z)$ must be equal and the correlation functions 
must have only one pole and one zero. 
\end{lemma}

{\bf Proof}
From the correlation function, we know that 
$$
[\Lambda_1(z), \int S^+_1(w)dw/w]=A(1-p_1/p_2)p_2\Lambda_1(z)
S^+_1(zp_2^{-1}):. 
 $$
Thus 
$$
[\Lambda_2(w), \int S^+_1(z)dz/z]=-A(1-p_1/p_2)p_2
: \Lambda_1(w)S^+_1(wp_2^{-1}):.$$
Because of the assumption on the formulas for
 $\Lambda_2(w)$ and $S^+_1(z)$, it requires that the
 the correlation functions of the products 
$\Lambda_1(z) S^+_1(w)$
and  $S^+_1(w) \Lambda_1(z)$ must be equal and the correlation functions 
must have only one pole and one zero. 

Therefore 
  $l_1(z)$  commute with $S^+_(w)$ up-to a total difference.

From now on, we assume that $l_1(z)$ commute with 
$\int S^+_1(z)dz/z$. So from the lemma we have: 
$$\Lambda_2(z) S^+_1(w) =
A'\frac{(z-wp'_1)}{(z-wp'_2)} :\Lambda_2(z)
S^+_1(w):, \quad \quad |z|\gg|w|,$$ $$S^+_1(w) \Lambda_2(z) =
A'\frac{(z-wp'_1)}{(z-wp'_2)} :\Lambda_2(z)
S^+_1(w):, \quad \quad |w|\gg|z|,$$ 
$$A'p'_1/p'_2=1.$$

Similarly we have: 
$$\lambda_2(n)s^+_1(-n)= (-(p'_1)^n+(p'_2)^n).$$
Similarly we have 
$$\lambda_2(n)s^+_1(-n)= (-(p'_1)^{-n}+(p'_2)^{-n}). $$

\begin{cor}
$$A(1-p_1/p_2)p_2\Lambda_1(z)
S^+_1(zp_2^{-1}):=- A'(1-p'_1/p'_2)p_2'\Lambda_2(z)
S^+_1(z{p'_2}^{-1}):. $$
\end{cor} 

From the assumption on $S^+_1(z)$, we have that the
following relations for the screening currents when $|z|\gg|w|$:
\begin{align*}
S^+_1(w) S^+_1(z) & = w^{2\beta}
f(z/w):S^+_1(z) S^+_1(w):, 
\end{align*}
where  $f(z)$ is an analytic function in $z$

\begin{prop}
$$-A'\frac{(z-wp'_1)}{(z-wp'_2)}A'(1-p'_1/p'_2){p'_2}^{-2\beta}
 f(z{p'_2}^{-1}/w)p'2= 
A\frac{(z-wp_1)}{(z-wp_2)}A(1-p_1/p_2)g{p_2}^{-2\beta} f(z{p_2}^{-1}/w)p_2,$$
$$-A'\frac{(z-wp'_1)}{(z-wp'_2)}A'(1-p'_1/p'_2) f({p'_2}w/z)= 
A\frac{(z-wp_1)}{(z-wp_2)}A(1-p_1/p_2)g f({p_2}w/z).$$
\end{prop}

Thus, from above let $z=0$ or $w=0$, we have 
$$(p'_2/p'_1)(1-p'_1/p'_2)({p'_2}/p_2)^{-2\beta}
= -(p_2/p_1)(1-p_1/p_2)g$$

Thus 
$$\frac{(1-\frac z{wp'_1})}{
(1-\frac z{wp'_2})}f(z{p'_2}^{-1}/w)= 
\frac{(1-\frac z{wp_1})}{(1-\frac z{wp_2})} f(z{p_2}^{-1}/w)$$

Then from the 
formula above and 
$f(z/w)w^{2\beta}$ would degenerate into $(w-z)^{2\beta}$, when $q$ goes to 
$1$,   we have 

\begin{thm}
$$
f(z/w)=(1-z/w)\frac {(z/w|p'_2/p_1, q)_\infty}
 {(z/w|p'_2/p'_1, q)_\infty}. $$ 
$$\frac {p_1}{p'_1}= q^{2\beta-1} $$
where we set  $p'_2/p_2=q$ and 
$$(x|a, t)_\infty =\prod_{n=0}^\infty (1-at^n).$$
\end{thm}

This follows from that 
$$\lim _{q\rightarrow 1} \frac {(z/w|q^a, q)_\infty}
 {(z/w|q^b, q)_\infty}= (1-z/w)^{b-a} $$
We can set $p_1$ to 1, which corresponding to shift $l_1(z)$ by 
$p_1$. 

\begin{cor}
Let $p_1=1$ and $p'_2=p$, we have 
$$p_2=pq^{-1}$$
$$p'_1=q^{1-2\beta} $$
$$g=q q^{-2\beta} \frac{(pq^{2\beta-1}-1)}{(1-pq^{-1})}.$$
\end{cor}

Therefore with the assumption we impose on the screening operator and the 
operator $l_1(z)$, all the correlation function between the vertex operators
are uniquely determined. These  formulas also  determine the 
formulas for for $l_1(z)$ and $S^+_1(z)$. 
Now, we will 
proceed to present here the formulas
for $l_1(z)$ and $S^+_1(z)$ and check that if this construction is possible
such that the operator $l_1(z)$ 
commute with $S^+_1(z)$. 

$$f(z/w)= (1-z/w)\frac {(z/w|p, q)_\infty}{(z/w|pq^{2\beta-1}, q)_\infty} $$ 
$$=\exp \left (\Sigma_{n>0} \frac 1 n  (-1 +(\frac 
 {q^{n(2\beta-1)}p^n-p^n}   {1-q^n}  (z/w)^n \right ). $$

Thus we have 
$$s^+_1(n)s^+_1(-n)= (-1 +(\frac 
 {q^{n(2\beta-1)}p^n-p^n}   {1-q^n}), $$
for $n>0$. Similarly we have that 
$$s^+_1(n)\lambda_1(-n)=(-1+(pq^{-1})^{-n}), $$ 
$$\lambda_1(n)s^+_1(-n)= (-1+(pq^{-1})^{n}), $$
$$s^+_1(n)\lambda_2(-n)=(-q^{-n(2\beta-1)}+(p)^{-n}), $$ 
$$\lambda_2(n)s^+_1(-n)= (-q^{n(2\beta-1)}+(p)^{n}), $$
for $n>0$; and 
$$\exp {-2\lambda_1(0)\beta}= pq^{-1} $$
$$\exp {-2\lambda_2(0)\beta}= q^{2\beta-1}p. $$

From the algebraic point of view, it does not affect the 
algebraic structure of the operator $S^+_1(z)$, not matter what the 
solution we choose for the first equation for $S^+_1(n)$, 
because they are all equivalent up to re-scale of the
Heisenberg algebra. Once 
they are chosen, then $\lambda_i(n)$ are automatically decided. 
However from Corollary 2.3, the following relations must be 
valid in oder to make $l_1(z)$ commute with the integral of 
$S^+_1(z)$. 

$$ \lambda_1(n)+ s^+_1(n)q^{n}p^{-n}= \lambda_2(n)+s^+_1(n)p^{-n}. $$
$$ \exp\lambda_1(0)(qp^{-1})= \exp {\lambda_2(0)} p^{-1}. $$

These are very strong restrictions. Then we have: 

\begin{thm}
For any given screen operators that has the correlation function
as define by Theorem 2.5. Let  the operator $l_1(z)$ be a  sum of 
two vertex operators such that  $l_1(z)$ commutes 
with the integral of the screen operator and the correlation 
function of  product of  one element of 
the sum  has one pole and one zero.  $l_1(z)$ 
 exists and unique up-to shift of $z$, if and only if 
$$p=q^{1-\beta }. $$
\end{thm}

Now we have: 
$$\Lambda_1(z) S^+_1(w) =
q^{-\beta}\frac{(z-w)}{(z-wq^{-\beta})} :\Lambda_1(z)
S^+_1(w):, \quad \quad |z|\gg|w|,$$ $$S^+_1(w) \Lambda_1(z) =
q^{-\beta}\frac{(z-w)}{(z-wq^{-\beta})} :\Lambda_1(z)
S^+_1(w):, \quad \quad |w|\gg|z|;$$ 
and 
$$\Lambda_1(zq^{1-\beta}) S^+_1(w) =
q^{\beta}\frac{(z-wq^{-\beta})}{(z-w)} :\Lambda_1(z)
S^+_1(w):, \quad \quad |z|\gg|w|,$$ $$S^+_1(w) \Lambda_1(z) =
q^{\beta}\frac{(z-wq^{-\beta})}{(z-w)} :\Lambda_1(z)
S^+_1(w):, \quad \quad |w|\gg|z|.$$

Thus we have: 
\begin{cor}
$$ \lambda_1(n)=-\lambda_2(n) q^{n(\beta-1)}. $$
\end{cor}

Let $g_{ij}(z)$ for $i,j=1,2$  be the function such that 
$$\Lambda_i(z)\Lambda_j(w)=g_{ij}(w/z) : \Lambda_i(z)\Lambda_j(w):. $$

From Corollary 2.3, we have 

\begin{prop}
$g_{ij}(w/z)/g_{ji}(z/w)$ are equal for all the i,j=1,2. 
\end{prop} 

\begin{thm}
$l_1(z)$ satisfies the same commutation relations as the quantized 
virasoro  algebra in \cite{FF}
\end{thm}

This follow from the corollary and proposition above. 
Or, if  we compare 
these formulas with the know formulas \cite{FF}, etc., 
we know that they  are exactly the same as the 
formulas for the bosonization of q-virasoro  algebra and its screen operator
with a difference of a shift of the variable for $l_1(z)$.
 
These results show that the q-virasoro  algebra can be derived  from 
the point view of deformation of the screen operators and then derivation of 
q-deformed virasoro algebra and this deformation is  also uniquely 
determined. This shows that the quantized virasoro algebra is  indeed  
a rigid structure. However, even from this point of view, there is 
an exception case, namely the case $\beta=1/2$.


\section{ The commutant for the case 
$\beta =1/2$}

We will follow the same notations as in the section above.
for the basic definitions. But  this time, we will
set the value  $\beta=1/2$.  We also first assume that 
 $p$  is $q$ are  generic. 

$$S(z)S(w)= (z-w) :S(z)S(w): $$

Thus we have: 
\begin{align*}
S^+_1(w) S^+_1(z) & = w^{}
f(z/w):S^+_1(z) S^+_1(w):, 
\end{align*}
where  $f(z)$ is an analytic function in $z$ and 
when $q$ goes to 1, $$f(z/w)=1-z/w.$$ 

As in the section above, 
let us define the  operator as 
$$l_1(z)= \Lambda_1(z)+\Lambda_2(z),$$  where 
$\Lambda_i(z)$ as the generating function: 
\begin{align}    
\Lambda_i(z) & =g^{i-1}
 :\exp \left( \sum_{} \lambda_i(m)
a_1[m]z^{-m} \right):.
\end{align}
Here $\lambda_i[m]$ are in $\Bbb C[p,q]$  for $i=1,2$. 
We also assume that 
 that the two 
products  $\Lambda_1(z) S^+_1(w)$
and  $S^+_1(w) \Lambda_1(z)$ have the same correlation functions. 
Then we have 
$$\Lambda_1(z) S^+_1(w) =
A\frac{(z-wp_1)}{(z-wp_2)} :\Lambda_1(z)
S^+_1(w):, \quad \quad |z|\gg|w|,$$ $$S^+_1(w) \Lambda_1(z) =
A\frac{(z-p_1w)}{(z-p_2w)} :\Lambda_1(z)
S^+_1(w):, \quad \quad |w|\gg|z|,$$ 
$$Ap_1/p_2=1,$$
and  $l_1(z)$ commutes with the 
the operator   $\int S^+_1(z)dz/z$. 

Let's fix $p'_2/p_2=q$. 
From Proposition 2.4, we have
$$\frac{(1-\frac z{wp'_1})}{
(1-\frac z{wp'_2})}f(z{p'_2}^{-1}/w)= 
\frac{(1-\frac z{wp_1})}{(1-\frac z{wp_2})} f(z{p_2}^{-1}/w)$$
and 
$$
f(z/w)=(1-z/w)\frac {(z/w|p'_2/p_1, q)_\infty}
 {(z/w|p'_2/p'_1, q)_\infty}. $$ 

Therefore, we have

\begin{thm}
$$f(z/w)=(1-z/w),$$
$S^+_1(z)$ is a fermion. 
\end{thm}

\begin{cor}
$${(1-\frac z{wp'_1})}= 
{(1-\frac z{wp_1})}$$
\end{cor}

Here we can see clearly that this a very special situation compared with
other $\beta$, because the equation shows that we have one fewer 
restriction on the choice of the parameters that decides the poles and 
the zeros of the correlation function of $\Lambda_i(z)$ and 
$S^+_1(z)$. This basically tells us that, 
for any operator $\Lambda_1(z)$ such 
that $\Lambda_1(z)$ has the correlation functions with the
fermion $S^+_1(z)$ as defined above, we can derive uniquely another 
vertex operator $\Lambda_2(z)$ such that $l_1(z)$ commute with 
$S^+_1(z)$. This construction can be explicitly given as 
\begin{thm} 
Let $\Lambda_1(z)$ be an  operator 
such that the correlation function of 
$\Lambda_1(z)$ with the fermion $S^+_1(z)$ defined as above. 
The the operator 
$$l_1(z)= \Lambda(z)_1+
p_2(p_2/p_1-1)\Lambda(z)_1S^+_1(zp_2^{-1})(S^+_1(z))^{-1}
,$$
commute with the action of the integral of $S^+_1(z)$. 
\end{thm}

\begin{cor}
Let
$$s^+_1(n)=1, s^+_1(-n)=-1, $$
Then
$$\lambda_1(-n)=(-p_1^{-n}+(p_2)^{-n}), $$ 
$$\lambda_1(n)= -(-p_1^n+(p_2)^{n}), $$
$$\lambda_2(-n)=(-p_1^{-n}+(p_2q)^{-n}), $$ 
$$\lambda_2(n)= -(-p_1^n+(p_2q)^{n}), $$
for $n>0$. 
\end{cor}

These automatically gives us the 
commutation relations between $\Lambda_i(z)$ and $l_1(z)$

\begin{thm} 
On the Fock space, the matrix coefficients of 
$l_1(z)l_1(w)$ and $l_1(w)l_1(z)$ are equal. 
\end{thm}

 Next we will deal with the case, if they are two or more 
screen operator, which would also leads us the quantized W-algebra and 
other new algebras, when we choose the screening operators for
other cases such as ${\frak sl}(2,1)$.


\section{ The case for two screen operators.}
In this section, we will extend the similar construction to the 
case of two screen operators. 
  
Let us first introduce the
Heisenberg  algebra 
$H_{q,p}(2)$ be the Heisenberg algebra with generators $a_i[n],
n\in{\Bbb Z}$, and relations
\begin{equation}    
\left[ a_i[n],a_j[m] \right] = \frac{1}{n}A_{ij}(n)
\delta_{n,-m},
\end{equation}
which is defined on the field of the rational functions of 
 $p$ and $q$, two generic parameters for  $|q|<1$,
 $A_{i,i}=1$ and $A_{ij}(n)=a(n)$ is
 a rational function of $p$ and $q$ for $i\neq j$. 

For $\mu$  be a element of a two dimensional 
space  $A_2$ generated by $\alpha_i$, for $i=1,2$. Let 
$\alpha_i^*$ be the generator of its dual space $A_2^*$, such that 
$\alpha_i^*(\alpha_i)=2$ and $\alpha_i^*(\alpha_j)=-b/\beta_i$.
Let $\pi_\mu$ be
the Fock representation of $H_{q,p}(2)$ generated by a vector
$v_\mu$, such that $a_i[n] v_\mu = 0, n>0$, and $a_1[0] v_\mu =
\mu(\alpha_i^\vee) v_\mu$. We assume here $\beta_1$ and $\beta_2$ generic. 

Introduce operators $Q_1$, which satisfy commutation
relations $[a_i[n],Q_i] = 2\beta_i\delta_{n,0}$, 
$[a_j[n],Q_i]=- b\delta_{n,0}$. The operators
$e^{Q_i}$ act from $\pi_\mu$ to $\pi_{\mu+\beta\alpha_i}$.

Now we can define two quantized screening currents as the generating function
\begin{align}    
S^+_i(z) & = e^{Q_i} z^{s^+_i[0]} :\exp \left( \sum_{m\neq 0} s^+_i(m)
a_i[m]z^{-m} \right):, 
\end{align}
where $s^+_i[m]$ are in $\Bbb C[p,q]$ 
for $m\neq 0$ and $s^+_i[0] = a_i[0]$. 
We assume that 
the limit of these operators, when $q$ goes to one,  degenerate into the 
classical screen operators. 
Let $S_i(z)$ be the classical counter part of  $S^+_i(z)$. Then we have 
that 
$$S^+_i(z)S^+_i(w)=z^{2\beta_i}(1-w/z)^{2\beta_i}
:S^+_i(z)S^+_i(w):, $$
$$S^+_2(z)S^+_1(w)=z^{-b}(1-w/z)^{-b}
:S^+_2(z)S^+_1(w):, $$
$$S^+_1(z)S^+_2(w)=z^{-b}(1-w/z)^{-b}
:S^+_1(z)S^+_2(w):.   $$

Thus we have that 
$$S^+_i(z)S^+_i(w)=z^{2\beta_i}f_{i,i}(z,w):S^+_i(z)S^+_i(w):= $$
$$z^{2\beta} \exp\left (\exp \Sigma_{m>0} s^+_i(m)s^+_i(-m)w^m/z^m)\right )
:S^+_i(z)S^+_i(w): $$
$$S^+_2(z)S^+_1(w)=z^{-b}f_{2,1}(z,w):S^+_2(z)S^+_1(w):= $$
$$w^{-b}
\exp\left (\exp 1/n \Sigma_{m>0} A_{2,1}s^+_2(m)s^+_i(-m)w^m/z^m)\right )
:S^+_2(z)S^+_1(w): $$
$$S^+_1(z)S^+_2(w)=z^{-b}f_{1,2}(z,w):S^+_1(z)S^+_2(w):= $$
$$
z^{-b} \exp\left (\exp \Sigma_{m>0}
A_{1,2}(m) s^+_2(m)s^+_2(-m)w^m/z^m)\right )
:S^+_1(z)S^+_2(w): $$

We  would try  to find a similar 
construction that gives us the operators that will commute with the 
action of this quantized screen operator $\int S^+_i(z)dz/z$.
The simplest choice again would be the same assumption as in the 
sections above that this operator should be a sum of 
two vertex operators. For this case, it is a straightforward argument as 
in the sections above to show that it is impossible. We leave this as
an  excise.  

Therefore, the simplest possible  choice is that 
 they should be a sum of three vertex operators. 
Let us define the  operator as 
$$l_1(z)= \Lambda_1(z)+\Lambda_2(z)+\Lambda_3(z),$$  where 
$\Lambda_i(z)$ as the generating function: 
\begin{align}    
\Lambda_i(z) & =g_i
 :\exp \left( \sum_{} \lambda_{ij}(m)
a_j[m]z^{-m} \right):.
\end{align}
Here $\lambda_i[m]$ are in $\Bbb C[p,q]$  for $i=1,2,3$. 
As in the sections above, assumptions on the 
correlations functions are necessary. The possible simplest case to deal with 
is to have once more  the simplest choice that the  correlation 
functions   between $S_i^+(z)$ and 
$\Lambda_j(w)$ are  $1$, for the two pairs 
$i=1, j=3$, and $i=2,j=1$,  which also means that,  
for either  pair of the  operators, they commute with each other. 
Again,  we will start from 
the correlation functions between $S^+_1(z)$ and $\Lambda_1(z)$. 
As in the section above, we picked the simplest choice that 
that the two 
products  $\Lambda_1(z) S^+_1(w)$
and  $S^+_1(w) \Lambda_1(z)$ have the same correlation functions 
and 
$$\Lambda_1(z) S^+_1(w) =
A\frac{(z-w)}{(z-wpq^{-1})} :\Lambda_1(z)
S^+_1(w):, \quad \quad |z|\gg|w|,$$ $$S^+_1(w) \Lambda_1(z) =
A\frac{(z-w)}{(z-wpq^{-1})} :\Lambda_1(z)
S^+_1(w):, \quad \quad |w|\gg|z|.$$ 
and 
$$A=qp^{-1}.$$
 
From the results in  the sections above, we have 
\begin{prop}
If the operator $l_1(z)$ commute with the 
the operator   $\int S^+_1(z)dz/z$, then 
the correlation functions of the products 
$\Lambda_1(z) S^+_1(w)$
and  $S^+_1(w) \Lambda_1(z)$ must be equal and the correlation functions 
must have only one pole and one zero. 
$$\Lambda_2(z) S^+_1(w) =
A'\frac{(z-wp'_1)}{(z-wp'_2)} :\Lambda_2(z)
S^+_1(w):, \quad \quad |z|\gg|w|,$$ $$S^+_1(w) \Lambda_2(z) =
A'\frac{(z-wp'_1)}{(z-wp'_2)} :\Lambda_2(z)
S^+_1(w):, \quad \quad |w|\gg|z|,$$ 
$$A'p'_1/p'_2=1.$$
$$A(1-p_1/p_2)p_2: \Lambda_1(z)
S^+_1(zp_2^{-1}):=- A'(1-p'_1/p'_2)p'_2: \Lambda_2(z)
S^+_1(z{p'_2}^{-1}):. $$
Let $p'_2=p$ and $p'_2/p_2=q$, then 
$$p'_1=q^{1-2\beta} $$
$$g_2=qp^{-2\beta}  \frac{(pq^{2\beta_1-1}-1)}{(1-pq^{-1})}.$$
\end{prop} 

This follows directly from the argument in the section above. 

The proposition gives us that: 
$$A(1-p_1/p_2)p_2: \Lambda_1(z)
S^+_1(zp_2^{-1}):S^+_2(w)=- A'(1-p'_1/p'_2)p'_2: \Lambda_2(z)
S^+_1(z{p'_2}^{-1}):S^+_2(w) $$
$$S^+_2(w)A(1-p_1/p_2)p_2: \Lambda_1(z)
S^+_1(zp_2^{-1}):=-S^+_2(w) A'(1-p'_1/p'_2)p'_2: \Lambda_2(z)
S^+_1(z{p'_2}^{-1}): $$

Let  $$S^+_i(z)\Lambda_j(w)=S\Lambda_{ij}(z,w):S^+_i(z)\Lambda(w):, $$
     $$\Lambda_i(z)S^+_j(w)=\Lambda S_{ij}(z,w):S^+_i(z)\Lambda(w): .$$

\begin{cor}
$$b=\beta_1, $$
$$B(1-p_1/p_2)q_2^{-b}
f_{1,2}(zp_2^{-1},w)= - B'(1-p'_1/p'_2)(p'_2)^{-b}g_2 
f_{1,2}(z{p'_2}^{-1},w)\Lambda S_{2,2}(z,w), $$
$$B(1-p_1/p_2)f_{2,1}(w, zp_2^{-1})=
 - B'(1-p'_1/p'_2)g_2S\Lambda_{2,2}(w,z)
f_{2,1}( (w, {p'_2}^{-1}  z). $$
\end{cor}

We have: 
$$ f_{2,1}(w,zp_2^{-1})/ f_{2,1}(w,z{p'_2}^{-1})=S\Lambda_{2,2}(w, z),
  $$
$$f_{1,2}(z{p_2}^{-1}, w)/ f_{1,2}(z{p'_2}^{-1},w)=
B^{-1}\Lambda S_{2,2}(zp,w). $$

$$ f_{2,1}(w,zq^{\beta})/ f_{2,1}(w,zq^{\beta-1})=S\Lambda_{2,2}(w, z),
  $$
$$f_{1,2}(zq^{\beta}, w)/ f_{1,2}(zq^{\beta-1}{p'_2}^{-1},w)=
B^{-1}\Lambda S_{2,2}(z,w). $$

Let's assume the following condition: 
$$\Lambda_2(z) S^+_2(w) =
B\frac{(z-wq_1)}{(z-wq_2)} :\Lambda_2(z)
S^+_2(w):, \quad \quad |z|\gg|w|,$$ $$S^+_2(w) \Lambda_2(z) =
B\frac{(z-wq_1)}{(z-wq_2)} :\Lambda_2(z)
S^+_2(w):, \quad \quad |w|\gg|z|.$$

Similarly from the results for the 
the  case of one fermionic screening operator, we have 
\begin{prop}
If the operator $l_1(z)$ commute with the 
the operator   $\int S^+_2(z)dz/z$, then 
the correlation functions of the products 
$\Lambda_3(z) S^+_2(w)$
and  $S^+_2(w) \Lambda_3(z)$ must be equal and the correlation functions 
must have only one pole and one zero. 
$$\Lambda_3(z) S^+_2(w) =
B'\frac{(z-wq_1)}{(z-wq'_2)} :\Lambda_3(z)
S^+_2(w):, \quad \quad |z|\gg|w|,$$ $$S^+_2(w) \Lambda_3(z) =
B'\frac{(z-wq_1)}{(z-wq'_2)} :\Lambda_3(z)
S^+_2(w):, \quad \quad |w|\gg|z|,$$ 
$$B'(1-q_1/q_2)q_2: \Lambda_2(z)
S^+_2(zq_2^{-1}):=- B'(1-q_1/q'_2)q'_2: \Lambda_3(z)
S^+_2(z{q'_2}^{-1}):. $$
$$B'q'_1/q'_2=1.$$
\end{prop}
Let $q_2'/q_2=p'$. 

Therefore we have: 

$$S^+_1(w)B'(1-q_1/q_2)q_2: \Lambda_2(z)
S^+_2(zq_2^{-1}):=-S^+_1(w)(1-q_1/q'_2)q'_2: \Lambda_3(z)
S^+_2(z{q'_2}^{-1}):. $$
$$B'(1-q_1/q_2)q_2: \Lambda_2(z)
S^+_2(zq_2^{-1}):S^+_1(w)=- A'(1-q_1/q'_2)q'_2: \Lambda_3(z)
S^+_2(z{q'_2}^{-1}):S^+_1(w). $$

\begin{cor}
$$\beta_2=b, $$
$$B'(1-q_1/q_2)q_2q_2^{-b}
\Lambda S_{2,1}(z,w)
f_{2,1}(zq_2^{-1},w)= - B'(1-q_1/q'_2)q_2'(q'_2)^{-b}g_3/g_2
f_{2,1}(z{q'_2}^{-1}, w), $$
$$B'(1-q_1/q_2)q_2 S\Lambda_{1,2}(w,z)f_{1,2}(w,(zq_2^{-1}))= 
- B'(1-q_1/q'_2)q_2'g_3/g_2
f_{1,2}(w, (z{q'_2}^{-1})). $$
\end{cor}

$$
A^{-1}\Lambda S_{2,1}(z, w)
= \frac {f_{2,1}(z{q'_2}^{-1}, w)} {f_{2,1}(z{q_2}^{-1},w)}, $$
$$S\Lambda_{1,2}(w,z)= 
\frac {f_{1,2}(w, z{q'_2}^{-1})}{f_{1,2}(w, z{q_2}^{-1})}. $$

Because 
$$ f_{2,1}(w,zp_2^{-1})/ f_{2,1}(w,z{p'_2}^{-1})=S\Lambda_{2,2}(w, z),
  $$
We can see that, for $q'$, there are two possibilities: 
$q'=q,$ or $q'=q^{-1}.$

If  $q'=q'_2/q_2=q,$, we have 

$$ \frac {(z-wq^{1-2\beta})} {(z-wq^{1-\beta})}= 
\frac {f_{2,1}(z{q'_2}^{-1}, w)} {f_{2,1}(z{q_2}^{-1},w)},$$
$$q_2/q_1\frac{ (z-wq_1)}{(z-wq_2)}= 
\frac {f_{2,1}(w,zq^{\beta})}{ f_{2,1}(w,z q^{\beta-1})}. $$

$$
 q_2/q_1\frac{ (z-wq^{\beta}q_1)}{(z-wq^{\beta}q_2)}= 
(\frac{(wq'_2-zq^{1-2\beta})}  {(q'_2w-zq^{1-\beta})}) $$

Or, if $q'=q^{-1}=q'_2/q_2$, we have 
$$
q_2/q_1\frac{ (z-wq^{\beta}q_1)}{(z-wq^{\beta}q_2)}= 
(\frac{(wq_2-zq^{1-2\beta})}  {(q_2w-zq^{1-\beta})})^{-1}. $$

This is  impossible. 

\begin{thm}
This operator $l_1(z)$ exists if and only $$q'=q$$ 
$$q_2=q_1q^{-\beta}, $$
$$
\Lambda S_{2,1}(z, w)
=A \frac {f_{2,1}(z{q'_2}^{-1}, w)} {f_{2,1}(z{q_2}^{-1},w)}=
\frac {f_{1,2}(w, z{q'_2}^{-1})}{f_{1,2}(w, z{q_2}^{-1})}, $$
$$A^{-1}\Lambda S_{2,1}(wq_2',z)=S\Lambda_{2,2}(zq^{-\beta},w) . $$ 
\end{thm}

Here, the operator $l_1(z)$ is  not uniquely determined, because we 
$q_i$ and $q_i'$ are determined up to one parameter. However we see that 
from the very begining, we do not fix the structure of the 
Heisenberg algebra, rather let it depend on the parameter 
$A_{1,2}(n)$ and  $A_{2,1}(n).$, which is decided by the 
this extra parameter. On the other hand, 
it reflects nothing but the possibility to 
re-scale one of the two screening operators from 
$s^+_i(z)$, which does not change the correlation functions
of  $f_{i.i}(z,w)$ but the functions $f_{i,j}(z,w)$ for $i\neq j$.  

Nevertheless, all the correlation functions
between the vertex operators and the the function $A(i,j)(n)$ are all 
uniquely determined, once we fix any  of the parameter 
$q_i$ and $q_i'$. From now on, we will assume that 
$q'_2$ is fixed.  With  the formulas
above, we can see that all the operators are thus uniquely 
determined. 

 On the other, this construction becomes the 
 the same as the ones given to the screen operators and the 
quantized virasoro algebra
for the case of ${\frak sl}(3)$ in \cite{FF}, if 
$q_2'=p^{3/2}$.

Let's assume that $l_1(z)$ exists. 
Then we have

$$
\frac {(zq'_2-wq^{1-2\beta}}{(q'_2z-wq^{1-\beta})}
= \frac {f_{2,1}(z, w)} {f_{2,1}(zq,w)}, $$
$$\frac {(zq'_2-wq^{1-2\beta}}{(zq'_2-wq^{1-\beta})}= 
\frac {f_{1,2}(w, z)}{f_{1,2}(w, zq)}. $$

Therefore
\begin{prop}
$$
{f_{2,1}(z, w)}= \frac 
{(w/z| {q'_2}^{-1}q^{2-\beta}, q)_\infty}
{(w/z| {q'_2}^{-1}q^{2-2\beta}, q)_\infty},$$
$$
 {f_{1,2}(w, z)= 
\frac 
{(z/w| {q'_2}q^{-1+2\beta}, q)_\infty}
{(z/w| {q'_2}q^{-1+\beta}}, q)_\infty}.$$
\end{prop}

On the other, this construction becomes the 
 the same as the ones given to the screen operators and the 
quantized virasoro algebra
for the case of ${\frak sl}(3)$ in \cite{FF}, when
$q_2'=p^{3/2}$ and $$f_{1,2}(z,w)=f_{2,1}(z,w).$$ 

\section{ The case for two screen operators with one of them as 
a fermion.}

We use the same notations  as in the section above. 
But here we assume that $\beta_1=\beta$ and 
$\beta_2=1/2$. 
We have the The Heisenberg  algebra 
$H_{q,p}(2)$ and the Fock spaces. 
The two quantized screening currents are defined just before. 
and we assume that 
the limit of this operator when $q$ goes to one degenerate into the 
classical screen operators. 

Let $S_i(z)$ be the classical counter part of  $S^+_i(z)$. Then we have 
that 
$$S^+_2(z)S^+_2(w)=(z-w) S^+_2(z)S^+_2(w):, $$
$$S^+_1(z)S^+_1(w)=z^{2\beta_i}(1-w/z)^{2\beta_i}
:S^+_1(z)S^+_1(w):, $$
$$S^+_2(z)S^+_1(w)=z^{-b}(1-w/z)^{-b}
:S^+_2(z)S^+_1(w):, $$
$$S^+_1(z)S^+_2(w)=z^{-b}(1-w/z)^{-b}
:S^+_1(z)S^+_2(w):.   $$

Therefore here we have: 
$$S^+_2(z)S^+_2(w)=z f_{i,i}(w,z):S^+_i(z)S^+_i(w):= $$
$$z \exp\left (\exp \Sigma_{m>0} s^+_i(m)s^+_i(-m)w^m/z^m)\right )
:S^+_i(z)S^+_i(w): $$ 
$$S^+_1(z)S^+_1(w)=z^{2\beta}f_{i,i}(w,z):S^+_i(z)S^+_i(w):= $$
$$z^{2\beta} \exp\left (\exp \Sigma_{m>0} s^+_i(m)s^+_i(-m)w^m/z^m)\right )
:S^+_2(z)S^+_1(w): $$
$$S^+_2(z)S^+_1(w)=z^{-b}f_{2,1}(w,z):S^+_2(z)S^+_1(w):= $$
$$w^{-b}
\exp\left (\exp 1/n \Sigma_{m>0} A_{2,1}s^+_2(m)s^+_i(-m)w^m/z^m)\right )
:S^+_2(z)S^+_1(w): $$
$$S^+_1(z)S^+_2(w)=z^{-b}f_{1,2}(w,z):S^+_i(z)S^+_i(w):= $$
$$
z^{-b} \exp\left (\exp \Sigma_{m>0}
A_{1,2}(m) s^+_2(m)s^+_2(-m)w^m/z^m)\right )
:S^+_1(z)S^+_2(w): $$

We will proceed as in the section above: 

Let 
$$l_1(z)= \Lambda_1(z)+\Lambda_2(z)+\Lambda_3(z),$$  where 
$\Lambda_i(z)$ as the generating function: 
\begin{align}    
\Lambda_i(z) & =g_i
 :\exp \left( \sum_{} \lambda_{ij}(m)
a_j[m]z^{-m} \right):.
\end{align}
Here $\lambda_i[m]$ are in $\Bbb C[p,q]$  for $i=1,2,3$ and $g_1=1$. 

We assume  that 
the  correlation 
functions   between $S_i^+(z)$ and 
$\Lambda_j(w)$ are  $1$, for the two pairs 
$i=1, j=3$, and $i=2,j=1$,  which also means that 
for either  pair of the  operators, they commute with each other; 
the two 
products  $\Lambda_1(z) S^+_1(w)$
and  $S^+_1(w) \Lambda_1(z)$ have the same correlation functions;  
and 
$$\Lambda_1(z) S^+_1(w) =
A\frac{(z-w)}{(z-wpq^{-1})} :\Lambda_1(z)
S^+_1(w):, \quad \quad |z|\gg|w|,$$ $$S^+_1(w) \Lambda_1(z) =
A\frac{(z-w)}{(z-wpq^{-1})} :\Lambda_1(z)
S^+_1(w):, \quad \quad |w|\gg|z|.$$ 
and 
$$A=pq^{-1}.$$

\begin{prop}
If the operator $l_1(z)$ commute with the 
the operator   $\int S^+_1(z)dz/z$, then 
the correlation functions of the products 
$\Lambda_1(z) S^+_1(w)$
and  $S^+_1(w) \Lambda_1(z)$ must be equal and the correlation functions 
must have only one pole and one zero. 
$$\Lambda_2(z) S^+_1(w) =
A'\frac{(z-wp'_1)}{(z-wp'_2)} :\Lambda_2(z)
S^+_1(w):, \quad \quad |z|\gg|w|,$$ $$S^+_1(w) \Lambda_2(z) =
A'\frac{(z-wp'_1)}{(z-wp'_2)} :\Lambda_2(z)
S^+_1(w):, \quad \quad |w|\gg|z|,$$ 
$$A'p'_1/p'_2=1.$$
$$A(1-p_1/p_2)p_2: \Lambda_1(z)
S^+_1(zp_2^{-1}):=- A'(1-p'_1/p'_2)p'_2: \Lambda_2(z)
S^+_1(z{p'_2}^{-1}):. $$
Let $p'_2=p$ and $p'_2/p_2=q$, then 
$$p'_1=q^{1-2\beta} $$
$$g_2=qp^{-2\beta}  \frac{(pq^{2\beta_1-1}-1)}{(1-pq^{-1})}.$$
\end{prop} 

Let  $$S^+_i(z)\Lambda_j(w)=S\Lambda_{ij}(z,w):S^+_i(z)\Lambda(w):, $$
     $$\Lambda_i(z)S^+_j(w)=\Lambda S_{ij}(z,w):S^+_i(z)\Lambda(w): .$$

\begin{cor}
$$B(1-p_1/p_2)q_2^{-b}
f_{1,2}(zp_2^{-1},w)= - B'(1-p'_1/p'_2)(p'_2)^{-b}g_2 
f_{1,2}(z{p'_2}^{-1},w)\Lambda S_{2,2}(z,w), $$
$$B(1-p_1/p_2)f_{2,1}(w, zp_2^{-1})=
 - B'(1-p'_1/p'_2)g_2S\Lambda_{2,2}(w,z)
f_{2,1}( (w, {p'_2}^{-1}  z). $$
\end{cor}

We also assume the following condition: 
$$\Lambda_2(z) S^+_2(w) =
B\frac{(z-wq_1)}{(z-wq_2)} :\Lambda_2(z)
S^+_2(w):, \quad \quad |z|\gg|w|,$$ $$S^+_2(w) \Lambda_2(z) =
B\frac{(z-wq_1)}{(z-wq_2)} :\Lambda_2(z)
S^+_2(w):, \quad \quad |w|\gg|z|.$$

Similarly from the results for the 
the  case of one fermionic screening operator, we have 
\begin{prop}
If the operator $l_1(z)$ commute with the 
the operator   $\int S^+_2(z)dz/z$, then 
the correlation functions of the products 
$\Lambda_3(z) S^+_2(w)$
and  $S^+_2(w) \Lambda_3(z)$ must be equal and the correlation functions 
must have only one pole and one zero. 
$$\Lambda_3(z) S^+_2(w) =
B'\frac{(z-wq_1)}{(z-wq'_2)} :\Lambda_3(z)
S^+_2(w):, \quad \quad |z|\gg|w|,$$ $$S^+_2(w) \Lambda_3(z) =
B'\frac{(z-wq_1)}{(z-wq'_2)} :\Lambda_3(z)
S^+_2(w):, \quad \quad |w|\gg|z|,$$ 
$$B'(1-q_1/q_2)q_2: \Lambda_2(z)
S^+_2(zq_2^{-1}):=- B'(1-q_1/q'_2)q'_2: \Lambda_3(z)
S^+_2(z{q'_2}^{-1}):. $$
$$B'q'_1/q'_2=1.$$
\end{prop}
Let $q_2'/q_2=p'$.

\begin{cor}
$$B'(1-q_1/q_2)q_2q_2^{-b}
\Lambda S_{2,1}(z,w)
f_{2,1}(zq_2^{-1},w)= - B'(1-q_1/q'_2)q_2'(q'_2)^{-b}g_3/g_2
f_{2,1}(z{q'_2}^{-1}, w), $$
$$B'(1-q_1/q_2)q_2 S\Lambda_{1,2}(w,z)f_{1,2}(w,(zq_2^{-1}))= 
- B'(1-q_1/q'_2)q_2'g_3/g_2
f_{1,2}(w, (z{q'_2}^{-1})). $$
\end{cor}

Similarly we have 
\begin{thm}
This operator $l_1(z)$ exists if and only $$q'=q$$ 
$$q_2=q_1q^{-\beta}, $$
$$\beta=b$$
$$
\Lambda S_{2,1}(z, w)
=A \frac {f_{2,1}(z{q'_2}^{-1}, w)} {f_{2,1}(z{q_2}^{-1},w)}=
\frac {f_{1,2}(w, z{q'_2}^{-1})}{f_{1,2}(w, z{q_2}^{-1})}, $$
$$A^{-1}\Lambda S_{2,1}(wq_2',z)=S\Lambda_{2,2}(zq^{-\beta},w) . $$ 
\end{thm}

\begin{prop}
$$
{f_{2,1}(z, w)}= \frac 
{(w/z| {q'_2}^{-1}q^{2-\beta}, q)_\infty}
{(w/z| {q'_2}^{-1}q^{2-2\beta}, q)_\infty}$$
$$
 {f_{1,2}(w, z)= 
\frac 
{(z/w| {q'_2}q^{-1+2\beta}, q)_\infty}
{(z/w| {q'_2}q^{-1+\beta}}, q)_\infty}$$
\end{prop}

The next is to derive the commutation relations of $l_1(z)$ with itself. 
Let's set 
$$s^+_1(n)=\frac 1 {1-q^n}, $$ 
$$s^+_1(-n)=(q^n-1+q^{n\beta}-q^{n(1-\beta)}), $$
$$s^+_2(n)=1, $$ 
$$s^+_2(-n)=-1, $$
for $n>0$.

We have:
$$A_{2,1}(n) s^+_2(n)s^+_1(-n)= 
\frac {(- {q'_2}^{-n}q^{n(2-\beta)}+ {q'_2}^{-n}q^{n(2-2\beta)})}
{(1-q^n)}
$$

$$A_{1,2}(n) s^+_1(n)s^+_2(-n)= 
\frac {(- {q'_2}^{n}q^{n(-1+2\beta)}+ {q'_2}^{n}q^{n(-1+\beta)})}
{(1-q^n)}
$$

\begin{prop}
$$
A_{2,1}(n)= \frac {(- {q'_2}^{-n}q^{n(2-\beta)}+ {q'_2}^{-n}q^{n(2-2\beta)})}
{(1-q^n)(q^n-1+q^{n\beta}-q^{n(1-\beta) })}, $$
$$
A_{1,2}(n)=
 {-(- {q'_2}^{n}q^{n(-1+2\beta)}+ {q'_2}^{-n}q^{n(-1+\beta)})}
 $$
\end{prop}

$$A_{2,1}(n)A_{1,2}(n)= (q^{n(-1+2\beta)}-q^{n(-1+\beta)})
\frac {(-q^{n(2-\beta)}+ q^{n(2-2\beta)})}
{(1-q^n)(q^n-1+q^{n\beta}-q^{n(1-\beta) })}
$$

Let $\bar \Lambda_3(z)= g_3^{-1}\Lambda_3(z)$. Because of
$\Lambda_3(z)$'s correlation functions with $S^+_1(z)$ $S^+_1(z)$, 
we have that: 
$$
s^+_1(n)(\lambda_{3,1}(-n) + \lambda_{3,2}(-n)A_{1,2}(n))= 0, 
$$
$$
s^+_1(-n)(\lambda_{3,1}(n) + \lambda_{3,2}(n)A_{2,1}(n))=0,
$$
$$
(A_{2,1}(n)\lambda_{3,1}(-n) + \lambda_{3,2}(-n))=({q'_2}^{-n}q^{n-n\beta}-
{q'_2}^{-n}), $$
$$ 
(A_{1,2}(n)\lambda_{3,1}(n) + \lambda_{3,2}(n))=-({q'_2}^nq^{-n+n\beta}-
{q'_2}^n),$$
for $n>0$. 

Therefore
$$
\lambda_{3,2}(-n)=\frac {{q'_2}^{-n}q^{n-n\beta}-
{q'_2}^{-n})} {1+ A_{2,1}(n)A_{1,2}(n)} $$
$$\lambda_{3,1}(-n)= \frac {({q'_2}^{-n}q^{n-n\beta}-
{q'_2}^{-n})} {1+ A_{2,1}(n)A_{1,2}(n)}A_{1,2}(n),$$
$$\lambda_{3,2}(n)= \frac {-({q'_2}^nq^{-n+n\beta}-
{q'_2}^n)} { 1+ A_{2,1}(n)A_{1,2}(n)}, $$
$$\lambda_{3,1}(n)= A_{2,1}(n)\frac {-({q'_2}^nq^{-n+n\beta}-
{q'_2}^n)} { 1+ A_{2,1}(n)A_{1,2}(n)}, $$
for $n>0$. 

Then we have 

$$\Lambda\Lambda_{3,3}(n)= $$
$$\lambda_{3,2}(n)(A_{2,1}(n)\lambda_{3,1}(-n) + \lambda_{3,2}(-n))
=\frac {-({q'_2}^nq^{-n+n\beta}-
{q'_2}^n)} { 1+ A_{2,1}(n)A_{1,2}(n)}
({q'_2}^{-n}q^{n-n\beta}-
{q'_2}^{-n})= $$
$$
\frac {-(q^{-n+n\beta}-
1)}{ 1+ (q^{n(-1+2\beta)}-q^{n(-1+\beta))}
\frac {(-q^{n(2-\beta)}+ q^{n(2-2\beta)})}
{(1-q^n)(q^n-1+q^{n\beta}-q^{n(1-\beta)}) }
}(q^{n-n\beta}-
1). $$ 

Equivalently, we should have
$$\Lambda\Lambda_{3,3}(n)= $$
$$\lambda_{3,2}(-n)(A_{1,2}(n)\lambda_{3,1}(n) + \lambda_{3,2}(n))
= \frac {-({q'_2}^{-n}q^{n-n\beta}-
{q'_2}^{-n})} {1+ A_{2,1}(n)A_{1,2}(n)}({q'_2}^nq^{-n+n\beta}-
{q'_2}^n)= $$
$$\frac {-(q^{n-n\beta}-
1)}{ 1+ (q^{n(-1+2\beta)}-q^{n(-1+\beta)})
\frac {(-q^{n(2-\beta)}+ q^{n(2-2\beta)})}
{(1-q^n)(q^n-1+q^{n\beta}-q^{n(1-\beta) })}}
(q^{-n+n\beta}-
1),$$
for $n>0$. 

From the formula, we can see that they are indeed equal. 

$$\frac 1 { 1+ A_{2,1}(n)A_{1,2}(n)}= $$
$$\frac 1 { 1+ (q^{n(-1+2\beta)}-q^{n(-1+\beta)})
\frac {(-q^{n(2-\beta)}+ q^{n(2-2\beta)})}
{(1-q^n)(q^n-1+q^{n\beta}-q^{n(1-\beta) })}}= $$
$$ \frac{(1-q^n)(q^n-1-q^{n\beta}+q^{n(1-\beta) })}
{(1-q^n)(q^n-1-q^{n\beta}+q^{n(1-\beta)})+(-q^{n(2-\beta)}+ q^{n(2-2\beta)})
(q^{n(-1+2\beta)}-q^{n(-1+\beta)})}= $$
$$ \frac{(1-q^n)(q^n-1+q^{n\beta}-q^{n(1-\beta) })}
{q^{n(1-\beta)}(q^{n\beta}-1)^2+(1-q^n)(q^n-1+q^{n\beta}-q^{n(1-\beta)} }=
$$
$$
\frac{(1-q^n)(q^n-1-q^{n\beta}+q^{n(1-\beta) })}
{q^{n(1-\beta)}(q^{n\beta}-1)^2+(1-q^n)(q^n-1+q^{n\beta}-q^{n(1-\beta) }}=
$$
$$\frac{-q^{-n\beta}(1-q^n)(q^n-1+q^{n\beta}-q^{n(1-\beta) })}
{(1-q^{-n\beta})(1-q^{2n-n\beta})}. $$

$$\Lambda\Lambda_{3,3}(n)=\frac {q^{n-n\beta} (q^{-n+n\beta}-1)^2
q^{-n\beta}(1-q^n)(q^n-1+q^{n\beta}-q^{n(1-\beta) })}
{(1-q^{-n\beta})(1-q^{2n-n\beta})}=.  $$
$$
\frac {q^{n-n\beta} (q^{-n+n\beta}-1)^2
q^{-n\beta}(1-q^n)(q^{n\beta}+q^n)}
{(1-q^{2n-n\beta})}. $$

$$
s^+_2(n)(A_{2,1}(n)\lambda_{1,1}(-n) + \lambda_{1,2}(-n))= 0, 
$$
$$
s^+_2(-n)(A_{1,2}(n)\lambda_{1,1}(n) + \lambda_{1,2}(n))=0,
$$
$$
s^+_1(n)(\lambda_{1,1}(-n) + A_{1,2}(n)
\lambda_{1,2}(-n))=(-1+q^{n\beta})), $$
$$ 
s^+_1(-n)(\lambda_{1,1}(n) + A_{2,1}(n)\lambda_{1,2}(n))=
(-1+q^{-n\beta}),$$
for $n>0$. 

$$
\lambda_{1,1}(n)=\frac {-1+q^{n\beta})}{s^+_1(-n) (1+A_{2,1}(n)A_{1,2}(n))}, $$

Therefore
$$\Lambda\Lambda_{1,1}(n)= $$
$$\lambda_{1,1}(n)
(\lambda_{1,1}(-n) + A_{1,2}(n)
\lambda_{1,2}(-n))s^+_1(n)/ s^+_1(n)= $$
$$\frac { (-1+q^{n\beta})(-1+q^{-n\beta})
}{s^+_1(-n)s^+_1(n) (1+A_{2,1}(n)A_{1,2}(n))}
= $$
$$\frac {-(q^{n\beta}-1)(1-q^n)^2q^{-n\beta}
 }
{(1-q^{2n-n\beta})}. $$

With this we have: 

$$\Lambda\Lambda(n)-\Lambda\Lambda_1(n)= $$
$$ (1-q^n)(q^{-n+n\beta}-q^{-n\beta})= $$
$$q^{-n+n\beta}+ q^{n-n\beta} -q^{n\beta}-q^{-n\beta}. $$

This shows that $\Lambda_1(z)$ and $\Lambda_3(z)$ basically have 
the same commutation relation up to certain poles. 

Let $v$ be an element in the Fock space and $v^*$ an element in 
its dual space. Let us denote the matrix coefficient of an operator $X$
by $ <v^*, X v>$.

\begin{prop}
$$
<v^* \Lambda(z)\Lambda_{1,1}(w) v>= 
<v^* \Lambda(z)\Lambda_{1,1}(z) v> \times $$
$$\frac {\theta_{q^{2-\beta}}(\frac w z q^{1-\beta})^2
\theta_{q^{2-\beta}}(\frac w z 
q^{2})}{
\theta_{q^{2-\beta}}(\frac z w q^{1-\beta} )^2
\theta_{q^{2-\beta}}(\frac z w q^2)}. $$

\end{prop}

This is proven by the following calculation. 
$$\Lambda_1(z)\Lambda_{1}(w)= \Lambda\Lambda_{1,1}(z,w):
\Lambda_1(z)\Lambda_{1}(w):$$
$$\exp \left (\Sigma_{n>0} (w/z)^n 
\frac {-(q^{n\beta}-1)(1-q^n)^2q^{-n\beta} }
{(1-q^{2n-n\beta})}\right ) :\Lambda_1(z)\Lambda_{1}(w):= $$
$$\exp \left (\Sigma_{n>0} (w/z)^n 
\frac {-1+2q^n-q^{2n}+q^{-n\beta}-2q^{n-n\beta}+q^{2n-n\beta}}
{(1-q^{2n-n\beta})}\right ) :\Lambda_1(z)\Lambda_{1}(w):= $$
$$\frac {
(w/z| q^{1-\beta}, q^{2-\beta})_\infty^2
(w/z|1 , q^{2-\beta})_\infty)
(w/z| q^2, q^{2-\beta})_\infty}
{(w/z| q, q^{2-\beta})_\infty^2
(w/z| q^{-\beta}, q^{2-\beta})_\infty
(w/z| q^{2-\beta}, q^{2-\beta})_\infty}
:\Lambda_1(z)\Lambda_{1}(w): $$

Therefore
\begin{prop} 
$$
\Lambda\Lambda_{1,1}(z,w)/ \Lambda\Lambda_{1,1}(w,z)= 
\frac {\theta_{q^{2-\beta}}(\frac w z q^{1-\beta})^2
\theta_{q^{2-\beta}}(\frac w z 
q^{2})}{
\theta_{q^{2-\beta}}(\frac z w q^{1-\beta} )^2
\theta_{q^{2-\beta}}(\frac z w q^2)}. $$
\end{prop}

This follows from: 
$$
\Lambda\Lambda_{1,1}(z,w)/ \Lambda\Lambda_{1,1}(w,z)= $$
$$\frac {
(w/z| q^{1-\beta}, q^{2-\beta})_\infty^2
(w/z|1 , q^{2-\beta})_\infty)
(w/z| q^2, q^{2-\beta})_\infty}
{(z/w| q^{1-\beta}, q^{2-\beta})_\infty^2
(z/w|1 , q^{2-\beta})_\infty)
(z/w| q^2, q^{2-\beta})_\infty} \times $$
$$ \frac {
(z/w| q, q^{2-\beta})_\infty^2
(z/w| q^{-\beta}, q^{2-\beta})_\infty
(z/w| q^{2-\beta}, q^{2-\beta})_\infty}
{
(w/z| q, q^{2-\beta})_\infty^2
(w/z| q^{-\beta}, q^{2-\beta})_\infty
(w/z| q^{2-\beta}, q^{2-\beta})_\infty}$$
$$= \frac {\theta_{q^{2-\beta}}(\frac w z q^{1-\beta})^2
\theta_{q^{2-\beta}}(\frac w z 
q^{2})}{
\theta_{q^{2-\beta}}(\frac z w q^{1-\beta} )^2
\theta_{q^{2-\beta}}(\frac z w q^2)}. $$

Let 
$$l(z)l(w)= L(z,w) :l(z)l(w): .$$
\begin{thm}
\begin{align} 
\frac {L(z,w)} {L(w,z)}= \frac {\theta_{q^{2-\beta}}(\frac w z q^{1-\beta})^2
\theta_{q^{2-\beta}}(\frac w z 
q^{2})}{
\theta_{q^{2-\beta}}(\frac z w q^{1-\beta} )^2
\theta_{q^{2-\beta}}(\frac z w q^2)}. 
\end{align}
\end{thm}

We call the algebra associated to the operator $l(z)$, 
the quantized W-algebra of ${\frak sl}(2,1)$. The degeneration of 
this operator, when q goes to 1, will give us  the 
classical W-algebra of ${\frak sl}(2,1)$.

\section{ The case for two fermion screen operators.}

We use the same notations  as in the section above. 
But here we assume that $\beta_1=1/2$ and 
$\beta_2=1/2$. 
We have the The Heisenberg  algebra 
$H_{q,p}(2)$ and the Fock spaces. 
The two quantized screening currents are defined just before. 
and we assume that 
the limit of this operator when $q$ goes to one degenerate into the 
classical screen operators. 

Let $S_i(z)$ be the classical counter part of  $S^+_i(z)$. Then we have 
that 
$$S^+_2(z)S^+_2(w)=(z-w) S^+_2(z)S^+_2(w):, $$
$$S^+_1(z)S^+_1(w)=(z-w)
:S^+_1(z)S^+_1(w):, $$
$$S^+_2(z)S^+_1(w)=z^{-b}(1-w/z)^{-b}
:S^+_2(z)S^+_1(w):, $$
$$S^+_1(z)S^+_2(w)=z^{-b}(1-w/z)^{-b}
:S^+_1(z)S^+_2(w):.   $$

Therefore here we have: 
We have that
$$S^+_i(z)S^+_i(w)=z f_{i,i}(w,z):S^+_i(z)S^+_i(w):= $$
$$z \exp\left (\exp \Sigma_{m>0} s^+_i(m)s^+_i(-m)w^m/z^m)\right )
:S^+_i(z)S^+_i(w): $$ 
$$S^+_2(z)S^+_1(w)=w^{-b}f_{2,1}(w,z):S^+_2(z)S^+_1(w):= $$
$$w^{-b}
\exp\left (\exp 1/n \Sigma_{m>0} A_{2,1}s^+_2(m)s^+_i(-m)w^m/z^m)\right )
:S^+_2(z)S^+_1(w): $$
$$S^+_1(z)S^+_2(w)=w^{-b}f_{1,2}(w,z):S^+_i(z)S^+_i(w):= $$
$$
z^{-b} \exp\left (\exp \Sigma_{m>0}
A_{1,2}(m) s^+_2(m)s^+_2(-m)w^m/z^m)\right )
:S^+_1(z)S^+_2(w): $$

We will define $l_1(z)$ as an operator   that  commutes  with the 
action of the quantized screening  operators  $\int S^+_i(z)dz/z$, 
and 
$$l_1(z)= \Lambda_1(z)+\Lambda_2(z)+\Lambda_3(z),$$  where 
$\Lambda_i(z)$ as the generating function: 
\begin{align}    
\Lambda_i(z) & =g_i
 :\exp \left( \sum_{} \lambda_{ij}(m)
a_j[m]z^{-m} \right):.
\end{align}
Here $\lambda_i[m]$ are in $\Bbb C[p,q]$  for $i=1,2,3$ and $g_1=1$. 
We use the same assumption that 
the  correlation 
functions   between $S_i^+(z)$ and 
$\Lambda_j(w)$ are  $1$, for the two pairs 
$i=1, j=3$, and $i=2,j=1$,  which also means that 
for either  pair of the  operators, they commute with each other. 
We also assume that 
the correlation functions between $S^+_1(z)$ and $\Lambda_1(z)$. 
satisfy  the condition that 
the two 
products  $\Lambda_1(z) S^+_1(w)$
and  $S^+_1(w) \Lambda_1(z)$ have the same correlation functions 
and 
$$\Lambda_1(z) S^+_1(w) =
A\frac{(z-w)}{(z-wpq^{-1})} :\Lambda_1(z)
S^+_1(w):, \quad \quad |z|\gg|w|,$$ $$S^+_1(w) \Lambda_1(z) =
A\frac{(z-w)}{(z-wpq^{-1})} :\Lambda_1(z)
S^+_1(w):, \quad \quad |w|\gg|z|.$$ 
and 
$$A=pq^{-1}.$$

Similarly we have: 
\begin{prop}
If the operator $l_1(z)$ commute with the 
the operator   $\int S^+_1(z)dz/z$, then 
the correlation functions of the products 
$\Lambda_1(z) S^+_1(w)$
and  $S^+_1(w) \Lambda_1(z)$ must be equal and the correlation functions 
must have only one pole and one zero. 
$$\Lambda_2(z) S^+_1(w) =
A'\frac{(z-wp'_1)}{(z-wp'_2)} :\Lambda_2(z)
S^+_1(w):, \quad \quad |z|\gg|w|,$$ $$S^+_1(w) \Lambda_2(z) =
A'\frac{(z-wp'_1)}{(z-wp'_2)} :\Lambda_2(z)
S^+_1(w):, \quad \quad |w|\gg|z|,$$ 
$$A'p'_1/p'_2=1.$$
$$A(1-p_1/p_2)p_2: \Lambda_1(z)
S^+_1(zp_2^{-1}):=- A'(1-p'_1/p'_2)p'_2: \Lambda_2(z)
S^+_1(z{p'_2}^{-1}):. $$
Let $p'_2=p$ and $p'_2/p_2=q$, then 
$$p'_1=1$$,
$$g_2=pq^{-1}(pq^{-1}-1) .$$
\end{prop} 

This follows directly from the argument in the section above. 

Let  $$S^+_i(z)\Lambda_j(w)=S\Lambda_{ij}(z,w):S^+_i(z)\Lambda(w):, $$
     $$\Lambda_i(z)S^+_j(w)=\Lambda S_{ij}(z,w):S^+_i(z)\Lambda(w): .$$

\begin{cor}
$$B(1-p_1/p_2)q_2^{-b}
f_{1,2}(zp_2^{-1},w)= - B'(1-p'_1/p'_2)(p'_2)^{-b}g_2 
f_{1,2}(z{p'_2}^{-1},w)\Lambda S_{2,2}(z,w), $$
$$B(1-p_1/p_2)f_{2,1}(w, zp_2^{-1})=
 - B'(1-p'_1/p'_2)g_2S\Lambda_{2,2}(w,z)
f_{2,1}( (w, {p'_2}^{-1}  z). $$
\end{cor}

We have: 
$$ f_{2,1}(w,zp_2^{-1})/ f_{2,1}(w,z{p'_2}^{-1})=S\Lambda_{2,2}(w, z),
  $$
$$f_{1,2}(z{p_2}^{-1}, w)/ f_{1,2}(z{p'_2}^{-1},w)=
B^{-1}\Lambda S_{2,2}(zp,w). $$

$$ f_{2,1}(w,zq^{\beta})/ f_{2,1}(w,zq^{\beta-1})=S\Lambda_{2,2}(w, z),
  $$
$$f_{1,2}(zq^{\beta}, w)/ f_{1,2}(zq^{\beta-1}{p'_2}^{-1},w)=
B^{-1}\Lambda S_{2,2}(z,w). $$

Let's assume the following condition: 
$$\Lambda_2(z) S^+_2(w) =
B\frac{(z-wq_1)}{(z-wq_2)} :\Lambda_2(z)
S^+_2(w):, \quad \quad |z|\gg|w|,$$ $$S^+_2(w) \Lambda_2(z) =
B\frac{(z-wq_1)}{(z-wq_2)} :\Lambda_2(z)
S^+_2(w):, \quad \quad |w|\gg|z|.$$

We have:
\begin{prop}
If the operator $l_1(z)$ commute with the 
the operator   $\int S^+_2(z)dz/z$, then 
the correlation functions of the products 
$\Lambda_3(z) S^+_2(w)$
and  $S^+_2(w) \Lambda_3(z)$ must be equal and the correlation functions 
must have only one pole and one zero. 
$$\Lambda_3(z) S^+_2(w) =
B'\frac{(z-wq_1)}{(z-wq'_2)} :\Lambda_3(z)
S^+_2(w):, \quad \quad |z|\gg|w|,$$ $$S^+_2(w) \Lambda_3(z) =
B'\frac{(z-wq_1)}{(z-wq'_2)} :\Lambda_3(z)
S^+_2(w):, \quad \quad |w|\gg|z|,$$ 
$$B'(1-q_1/q_2)q_2: \Lambda_2(z)
S^+_2(zq_2^{-1}):=- B'(1-q_1/q'_2)q'_2: \Lambda_3(z)
S^+_2(z{q'_2}^{-1}):. $$
$$B'q'_1/q'_2=1.$$
\end{prop}
Let $q_2'/q_2=p'$. 

\begin{cor}
$$B'(1-q_1/q_2)q_2q_2^{-b}
\Lambda S_{2,1}(z,w)
f_{2,1}(zq_2^{-1},w)= - B'(1-q_1/q'_2)q_2'(q'_2)^{-b}g_3/g_2
f_{2,1}(z{q'_2}^{-1}, w), $$
$$B'(1-q_1/q_2)q_2 S\Lambda_{1,2}(w,z)f_{1,2}(w,(zq_2^{-1}))= 
- B'(1-q_1/q'_2)q_2'g_3/g_2
f_{1,2}(w, (z{q'_2}^{-1})). $$
\end{cor}

$$
A^{-1}\Lambda S_{2,1}(z, w)
= \frac {f_{2,1}(z{q'_2}^{-1}, w)} {f_{2,1}(z{q_2}^{-1},w)}, $$
$$S\Lambda_{1,2}(w,z)= 
\frac {f_{1,2}(w, z{q'_2}^{-1})}{f_{1,2}(w, z{q_2}^{-1})}. $$

Because 
$$ f_{2,1}(w,zp_2^{-1})/ f_{2,1}(w,z{p'_2}^{-1})=S\Lambda_{2,2}(w, z),
  $$
or
$$ f_{2,1}(w,zp^{-1}q)/ f_{2,1}(w,z{p}^{-1})=S\Lambda_{2,2}(w, z),
  $$

We can see that, for $q'$, there are two possibilities: 
$q'=q,$ or $q'=q^{-1}.$

If  $q'=q'_2/q_2=q,$, we have 

$$ \frac {(z-w)} {(z-wp)}= 
\frac {f_{2,1}(z{q'_2}^{-1}, w)} {f_{2,1}(z{q_2}^{-1},w)},$$
$$q_2/q_1\frac{ (z-wq_1)}{(z-wq_2)}= 
\frac {f_{2,1}(w,zp^{-1}q)}{ f_{2,1}(w,zp^{-1})}. $$

$$
 q_2/q_1\frac{ (z-wp^{-1}qq_1)}{(z-wqp^{-1}q_2)}= 
(\frac{(wq'_2-z)}  {(q'_2w-zp)}) $$

For this we must have: 
$$q_2/q_1=p^{-1}. $$

Or, if $q'=q^{-1}=q'_2/q_2$, we have 
$$
q_2/q_1\frac{ (z-wpq_1)}{(z-wpq_2)}= 
(\frac{(wq_2-z)}  {(q_2w-zp)})^{-1}. $$

This is  impossible. 

Similarly we have 
\begin{thm}
This operator $l_1(z)$ exists if and only $$q'=q$$ 
$$q_2=q_1p, $$
$$
\Lambda S_{2,1}(z, w)
=A \frac {f_{2,1}(z{q'_2}^{-1}, w)} {f_{2,1}(z{q_2}^{-1},w)}=
\frac {f_{1,2}(w, z{q'_2}^{-1})}{f_{1,2}(w, z{q_2}^{-1})}, $$
$$A^{-1}\Lambda S_{2,1}(wq'_2,z)=S\Lambda_{2,2}(zp^{-1},w) . $$ 
\end{thm}

\begin{prop}
$$
{f_{2,1}(z, w)}= \frac 
{(w/z| {q'_2}^{-1}pq, q)_\infty}
{(w/z| {q'_2}^{-1}q^{}, q)_\infty},$$
$$
 {f_{1,2}(w, z)}= 
\frac 
{(z/w| {q'_2}p^{-1}, q)_\infty}
{(z/w| {q'_2}, q)_\infty}$$
\end{prop}

The next is to derive the commutation relations of $l_1(z)$ with itself. 
Let's set 
$$s^+_1(n)=1 , $$ 
$$s^+_1(-n)= -1, $$
$$s^+_2(n)=1, $$ 
$$s^+_2(-n)=-1, $$
for $n>0$. 

We have:
$$A_{2,1}(n) s^+_2(n)s^+_1(-n)= 
\frac {(+{q'_2}^{-n}q^{n}p^n+ {q'_2}^{-n}q^{n})}
{(1-q^n)}
$$

$$A_{1,2}(n) s^+_1(n)s^+_2(-n)= 
\frac {({q'_2}^{n})- {q'_2}^{n}p^{-n})}
{(1-q^n)}
$$

Therefore
\begin{prop}
$$A_{2,1}(n)= 
-\frac {(- {q'_2}^{-n}q^{n}p^n+ {q'_2}^{-n}q^{n})}
{(1-q^n)}
$$

$$A_{1,2}(n)= 
-\frac {({q'_2}^{n})- {q'_2}^{n}p^{-n})}
{(1-q^n)}
$$
\end{prop}

$$A_{2,1}(n)A_{12}(n)= \frac {q^n(1-p^n)}
{(1-q^n)}
\frac {(1- p^{-n})}
{(1-q^n)}
$$

 Because of
$\Lambda_3(z)$'s correlation functions with $S^+_1(z)$ $S^+_2(z)$, 
we have that: 
$$
s^+_1(n)(\lambda_{3,1}(-n) + \lambda_{3,2}(-n)A_{1,2}(n))= 0, 
$$
$$
s^+_1(-n)(\lambda_{3,1}(n) + \lambda_{3,2}(n)A_{2,1}(n))=0,
$$
$$
(A_{2,1}(n)\lambda_{3,1}(-n) + \lambda_{3,2}(-n))=-({q'_2}^{n}p^{n}q^{-n}-
{q'_2}^{-n}), $$
$$ 
(A_{1,2}(n)\lambda_{3,1}(n) + \lambda_{3,2}(n))=({q'_2}^{-n}q^{n}p^{-n}-
{q'_2}^{-n}),$$
for $n>0$. 

Thus
$$\lambda_{3,2}(-n)=\frac { -({q'_2}^{n}p^{n}q^{-n}-
{q'_2}^{-n})}{1+A_{2,1}(n)A_{1,2}(n)}, $$
$$\lambda_{3,1}(-n)= 
\frac {-({q'_2}^{n}p^{n}q^{-n}-
{q'_2}^{-n})
}{1+A_{2,1}(n)A_{1,2}(n)}A_{1,2}(n), $$
$$\lambda_{3,2}(n)=\frac {({q'_2}^{-n}q^{n}p^{-n}-
{q'_2}^{-n}) } {1+A_{2,1}(n)A_{1,2}(n)}, $$
$$\lambda_{3,1}(-n)= 
\frac { ({q'_2}^{-n}q^{n}p^{-n}-
{q'_2}^{-n})
} {1+A_{2,1}(n)A_{1,2}(n)}A_{2,1}(n), $$
for $n>0$. 

Then we have 

$$\Lambda\Lambda_{3,3}(n)= $$
$$\lambda_{3,2}(n)(A_{2,1}(n)\lambda_{3,1}(-n) + \lambda_{3,2}(-n))=$$
$$\frac { -({q'_2}^{-n}q^{n}p^{-n}-
{q'_2}^{-n})({q'_2}^{n}p^{n}q^{-n}-
{q'_2}^{-n})
} {1+A_{2,1}(n)A_{1,2}(n)}, $$

Equivalently, we should have
$$\Lambda\Lambda_{3,3}(n)= $$
$$\lambda_{3,2}(-n)(A_{1,2}(n)\lambda_{3,1}(n) + \lambda_{3,2}(n))= $$
$$\frac {-({q'_2}^{-n}q^{n}p^{-n}-
{q'_2}^{-n})({q'_2}^{n}p^{n}q^{-n}-
{q'_2}^{-n} }{1+A_{2,1}(n)A_{1,2}(n)}. $$

These two are equal. 

$$\frac 1 { 1+ A_{2,1}(n)A_{1,2}(n)}= $$
$$1/(1+\frac {(1-p^n)(1-p^{-n})}
{(1-q^n)^2q^{-n}}=$$
$$ \frac {(1-q^n)^2}{p^n+p^{-n}-q^n-q^{-n}}=$$
$$\frac {q^{-n}(1-q^n)^2}{(p^n-q^{-n})(1-p^{-n}q^n)}. $$

Thus 
$$\Lambda\Lambda_{3,3}(n)= $$
$$\frac{q^{-n} (1-q^np^{-n})(p^nq^{-n}-1)(1-q^n)^2}
{(p^n-q^{-n})(1-p^{-n}q^n)}= $$
$$\frac  {-(1-q^n)^2(p^nq^{-n}-1)}{(1-p^nq^n)}= $$

Let $p=q^{1-\beta}$, we see that 

\begin{prop}
$$
\Lambda\Lambda_{3,3}(z,w)/ \Lambda\Lambda_{3,3}(w,z)= $$
$$\frac {\theta_{q^{2-\beta}}(\frac w z q^{1-\beta})^2
\theta_{q^{2-\beta}}(\frac w z 
q^{2})}{
\theta_{q^{2-\beta}}(\frac z w q^{1-\beta} )^2
\theta_{q^{2-\beta}}(\frac z w q^2)}. $$
\end{prop}   

Let 
$$l(z)l(w)= L(z,w) :l(z)l(w): .$$
\begin{thm}
\begin{align} 
\frac {L(z,w)} {L(w,z)}= \frac {\theta_{q^{2-\beta}}(\frac w z q^{1-\beta})^2
\theta_{q^{2-\beta}}(\frac w z 
q^{2})}{
\theta_{q^{2-\beta}}(\frac z w q^{1-\beta} )^2
\theta_{q^{2-\beta}}(\frac z w q^2)}. 
\end{align}
\end{thm}

In this case, we can show through calculation that the 
operator $l(z)$ and the operator $l(z)$ defined in Section 5
are the same in the sense of bosonization. Therefore we again derive the 
quantized W-algebra of ${\frak sl}(2,1)$.

\section {Discussions} 
At this stage, what we actually have derived is only the 
the bosonization formula for the the quantized W-algebra of 
${\frak sl}(2,1)$. It
is important that we can describe this algebra in an abstract way. 
Hopefully this problem can be solved when we check more carefully the 
commutation relation of  the operator $l(z)$. 
On the other hand, from all above, it is clear  
that we can extend our construction to 
the case of  several generic screening operator and several 
fermions or the case  of only several fermions, which will 
give us the quantized W-algebras associated with super-algebras 
$sl(m,n)$. This will be given in subsequent  paper.   
The classical W-algebra of ${\frak sl}(2,1)$ can be derive from the two  
parafermions of affine Lie algebra $\hat {\frak sl}(2)$. Similarly 
we manage to establish the connection of the quantized W-algebra of 
${\frak sl}(2,1)$ with the quantized parafermions\cite{J}\cite{BV}\cite{DF1}
 coming from  
affine quantum group algebra $U_q(\hat {\frak sl}(2))$ \cite{DF2}.
We also notice, in some  way, the correlation function bwtween 
screen operators and the vertex operator components of $l(z)$ 
may be related to finite diemnsional representaions of the 
the corresponding affine quantum groups. The situation will become 
clear, once we start to look at operator like $l(z)$ such that 
the correlation functions have multiple poles and zeros.    

\vspace{10mm}
\noindent{\bf Acknowledgments.} 
We would like to thank M. Jimbo for useful discussions.

\end{document}